\documentclass[a4paper]{amsart}
\usepackage{amssymb}
\usepackage{xypic}
\usepackage[only,mapsfrom]{stmaryrd}
\usepackage{graphicx}

\newtheorem{theorem}{Theorem}[section]
\newtheorem{lemma}[theorem]{Lemma}

\theoremstyle{definition}
\newtheorem{definition}[theorem]{Definition}
\newtheorem{example}[theorem]{Example}

\theoremstyle{remark}

\numberwithin{equation}{section}

\newcommand{\C}{\mathbb{C}}

\DeclareMathOperator{\Maps}{Maps}

\DeclareMathOperator{\link}{link}

\title[Positive scalar curvature]{Every quasitoric manifold admits an invariant metric of positive scalar curvature}

\author{Michael Wiemeler}
\address{Max-Planck-Institute for Mathematics, Vivatsgasse 7, D-53111 Bonn, Germany}
\email{wiemeler@mpim-bonn.mpg.de}
\thanks{The research was supported by a grant form the MPG}



\begin{document}
\begin{abstract}
  We prove that every quasitoric manifold admits an invariant metric of positive scalar curvature.
\end{abstract}

\maketitle


\section{Introduction}
\label{sec:intro}

A (strongly) quasitoric manifold is a \(2n\)-dimensional simply-connected manifold \(M\) which admits a smooth action of an \(n\)-dimensional compact torus \(T^n\) such that
\begin{itemize}
\item The \(T^n\)-action on \(M\) is locally standard, i.e. the \(T^n\)-action on \(M\) is locally modelled on the standard action of \(T^n\) on \(\C^n\).
\item If the first property is satisfied, then \(M/T\) is naturally an \(n\)-dimensional manifold with corners. We require that \(M/T\) is diffeomorphic to a simple polytope \(P\). In this case we call \(M\) a (strongly) quasitoric manifold over \(P\).
\end{itemize}

Quasitoric manifolds were first studied by Davis and Januszkiewicz \cite{0733.52006}.
An example of a quasitoric manifold is \(\C P^n\).
Moreover, every compact projective non-singular toric variety is a quasitoric manifold. 

The purpose of this note is to prove the following theorem.
\begin{theorem}
\label{sec:introduction}
  Every strongly quasitoric manifold of dimension greater than zero admits an \(T^n\)-invariant metric of positive scalar curvature.
\end{theorem}

Before we prove this theorem, we want to summarize what is known about the curvature of metrics on quasitoric manifolds in low dimensions.

In dimension two there is only one quasitoric manifold, namely \(\C P^1\). 
It admits an \(T^1\)-invariant metric of positive curvature.
By a result of Bazajkin and Matvienko \cite{1164.53349} every four-dimensional quasitoric manifold admits an invariant metric of positive Ricci-curvature.
In higher dimensions we will first use equivariant surgeries on moment-angle complexes to construct invariant metrics of positive scalar curvature on these manifolds.
Then Theorem~\ref{sec:introduction} follows from a result of Berard Bergery \cite{berard83:_scalar}.

This note is organized as follows.
In section \ref{sec:construct} we describe how equivariant surgery can be used to construct invariant metrics of positive scalar curvature.
In section \ref{sec:quasi_toric} we recall some basic facts about quasitoric manifolds and moment-angle complexes.
In section \ref{sec:bordism} we discuss equivariant surgeries of moment-angle complexes.
In section \ref{sec:proof} we prove Theorem \ref{sec:introduction}.

\section{Constructing metrics of positive scalar curvature}
\label{sec:construct}

Let \(G\) be a compact connected Lie-group and \(M\) a \(G\)-manifold.
In this section we recall from \cite{berard83:_scalar} and \cite{MR2376283} how equivariant surgery can be used to construct a \(G\)-invariant metric of positive scalar curvature on \(M\).
Note that if \(G\) is non-abelian and the \(G\)-action on \(M\) is effective, then there is a direct construction of an invariant metric of positive scalar curvature on \(M\) \cite{MR0358841}.

Now let \(N\) be a \(G\)-manifold of the same dimension as \(M\).
Let \(H\) be a closed subgroup of \(G\).
We say that \(M\) is obtained from \(N\) by an equivariant surgery of codimension \(c\), if there are orthogonal \(H\)-representations \(V_1\) and \(V_2\) with \(\dim V_1= \dim N-\dim G/H -c+1\), \(\dim V_2=c\)  and an equivariant embedding 
\begin{equation*}
  \phi:G\times_{H}(S(V_1)\times D(V_2))\rightarrow N,
\end{equation*}
such that \(M\) is equivariantly diffeomorphic to
\begin{equation*}
  (N-\phi(G\times_H(S(V_1)\times D(V_2))))\cup_{G\times_H(S(V_1)\times S(V_2))} G \times_H(D(V_1)\times S(V_2)).
\end{equation*}
Here \(S(V_i)\) and \(D(V_i)\) denote the unit sphere and the unit disc in \(V_i\), respectively.

Then there is the following equivariant version of a theorem which has been proven independently by Gromov and Lawson \cite{MR577131} and Schoen and Yau \cite{MR535700}.

\begin{theorem}[{\cite[Theorem 11.1]{berard83:_scalar}, \cite[Theorem 2]{MR2376283}}]
\label{sec:constr-metr-posit}
  Let \(M\) and \(N\) be \(G\)-manifolds.
  Assume that \(N\) admits an \(G\)-invariant metric of positive scalar curvature.
  If \(M\) is obtained from \(N\) by equivariant surgery of codimension at least three, then \(M\) admits an invariant metric of positive scalar curvature.
\end{theorem}

\section{Quasitoric manifolds and moment-angle complexes}
\label{sec:quasi_toric}

In this section we recall some basic facts about quasitoric manifolds and moment-angle complexes.
For more details we refer the reader to \cite{0733.52006} and \cite{MR1897064}.

Let \(P\) be an \(n\)-dimensional simple polytope.
We call the codimension-one faces of \(P\) facets of \(P\).
Denote by \(\mathfrak{F}(P)\) the set of facets of \(P\).
For each \(F_0\in \mathfrak{F}(P)\) denote by \(T_{F_0}\) the subtorus 
\begin{equation*}
  \{t\in T^{\mathfrak{F}(P)};\; t(F)=1 \text{ if } F\in\mathfrak{F}(P) \text{ and }F\neq F_0\}
\end{equation*}
of \(T^{\mathfrak{F}(P)}=\Maps(\mathfrak{F}(P), S^1)\).

Then the moment-angle complex of \(P\) is defined to be
\begin{equation*}
  Z_P=P\times T^{\mathfrak{F}(P)}/\sim,
\end{equation*}
where \((x,t)\sim(x',t')\) if and only if \(x=x'\) and \(tt'^{-1}\in \prod_{x\in F} T_F\).
The torus \(T^{\mathfrak{F}(P)}\) acts on \(Z_P\) by left translation on the second factor.

Up to equivariant diffeomorphism there is exactly one smooth structure on \(Z_P\) for which this \(T^{\mathfrak{F}(P)}\)-action is smooth and \(Z_P/T^{\mathfrak{F}(P)}\) is diffeomorphic to \(P\). (This follows from Corollary 4.7 of \cite{MR2285318} or by an argument which is similar to the one given in Section 5 of \cite{wiemeler12:_exotic} in the case of quasitoric manifolds.)
We equip \(Z_P\) with this smooth structure.

If \(M\) is a (strongly) quasitoric manifold over \(P\),
then there is a subtorus \(T'\) of \(T^{\mathfrak{F}(P)}\) such that \(T'\) acts freely on \(Z_P\) and \(M\) is diffeomorphic to \(Z_P/T'\).
Moreover, the torus action on \(M\) is given by the natural action of \(T^{\mathfrak{F}(P)}/T'\) on \(Z_P/T'\) \cite[Proposition 6.5]{MR1897064}.
\begin{example}
  If \(P=\Delta^n\) is an \(n\)-dimensional simplex, then \(T^{\mathfrak{F}(\Delta^n)}=T^{n+1}\) is an \(n+1\)-dimensional torus.
Moreover, \(Z_{\Delta^n}\) is equivariantly diffeomorphic to \(S^{2n+1}\), where the \(T^{n+1}\)-action on \(S^{2n+1}\) is induced by the standard action of \(T^{n+1}\) on \(\C^{n+1}\).
Denote by \(T'\) the diagonal subtorus of \(T^{n+1}\).
Then \(\C P^n=S^{2n+1}/T'\) with the induced action of \(T^{n+1}/T'\) is a quasitoric manifold over \(\Delta^{n}\). 
\end{example}

\section{Equivariant surgeries of moment-angle complexes}
\label{sec:bordism}

In this section we recall some facts about equivariant surgeries of moment-angle complexes.

If \(P\) is a simple polytope, then the dual polytope \(P^*\) of \(P\) is simplicial.
We denote by \(K(P)\) the boundary complex of \(P^*\).
It can be viewed as a simplicial complex on the vertex set \(\mathfrak{F}(P)\) such that \(\sigma\subset \mathfrak{F}(P)\) is a simplex of \(K(P)\) if and only if \(\bigcap_{F\in\sigma}F\neq \emptyset\).

For the proof of our main theorem we will need certain constructions on simplicial complexes so called bistellar moves.

\begin{definition}
  Let \(K\) be a pure \((n-1)\)-dimensional simplicial complex and \(\sigma\in K\) a \((n-1-i)\)-dimensional simplex, \(0\leq i\leq n-1\), such that \(\link_K \sigma\) is the boundary \(\partial\tau\) of a \(i\)-dimensional simplex that is not a face of \(K\).
  Then the operation \(\chi_\sigma\) defined by
  \begin{equation*}
    \chi_\sigma(K)=(K-(\sigma * \partial \tau))\cup (\partial \sigma * \tau)
  \end{equation*}
is called a bistellar \(i\)-move.
Here \(\sigma *\tau\) denotes the join of the two complexes \(\sigma\) and \(\tau\).
\end{definition}

There are the following relations between bistellar moves and equivariant surgeries of moment-angle complexes.

\begin{lemma}[{\cite[Example 6.22 and Construction 6.23]{MR1897064}}]
\label{sec:equiv-surg-moment-1}
  Let \(P\) and \(P'\) be simple \(n\)-dimensional convex polytopes of dimension \(n\geq 2\) and \(0\leq i\leq n-1\).
  Assume that \(K(P')\) can be obtained from \(K(P)\) by a bistellar \(i\)-move.
  Then we have:
  \begin{enumerate}
  \item If \(i=0\), then there is an isomorphism \(T^{\mathfrak{F}(P')}\rightarrow T^{\mathfrak{F}(P)}\times S^1\) and \(Z_{P'}\) can be obtained from \(Z_P\times S^1\) by equivariant surgery of codimension \(2n\).
  \item If \(1\leq i\leq n-2\), then there is an isomorphism \(T^{\mathfrak{F}(P')}\rightarrow T^{\mathfrak{F}(P)}\) and \(Z_{P'}\) can be obtained from \(Z_P\) by equivariant surgery of codimension \(2n-2i\).
  \item If \(i=n-1\), then there is an isomorphism \(T^{\mathfrak{F}(P')}\times S^1\rightarrow T^{\mathfrak{F}(P)}\) and \(Z_{P'}\times S^1\) can be obtained from \(Z_P\) by equivariant surgery of codimension \(2\).
\end{enumerate}
\end{lemma}

For the proof of Theorem \ref{sec:introduction} we need the following result of Ewald.

\begin{theorem}[\cite{MR510150}]
\label{sec:equiv-surg-moment}
  Let \(P\) be a simple polytope of dimension \(n\geq 3\).
  Then there is a sequence of simple polytopes \(P_1,\dots,P_m\) such that \(P_1=\Delta^n\) and \(P_m=P\) and, for \(i=1,\dots, m-1\), \(K(P_{i+1})\) is obtained from \(K(P_i)\) by a bistellar \(k\)-move with \(0\leq k\leq n-2\).
\end{theorem}

\section{Proof of Theorem \ref{sec:introduction}}
\label{sec:proof}

In this section we prove Theorem \ref{sec:introduction}.
Let \(P\) be a simple polytope of dimension \(n\) and \(M\) a quasitoric manifold over \(P\).
By the results cited in the introduction we may assume that \(n\geq 3\).
By Lemma~\ref{sec:equiv-surg-moment-1} and Theorem \ref{sec:equiv-surg-moment} we know that \(Z_P\) is obtained from \(Z_{\Delta^n}\times \prod S^1=S^{2n+1}\times \prod S^1\) by a sequence of equivariant surgeries of codimension greater or equal to four.
Since \(S^{2n+1}\times \prod S^1\) admits an obvious invariant metric of positive scalar curvature, it follows from Theorem \ref{sec:constr-metr-posit} that this also holds for \(Z_P\).

Now \(M\) is the quotient of a free action of a subtorus \(T'\) of  \(T^{\mathfrak{F}(P)}\) on \(Z_{P}\).
Hence, it follows from Theorem C of \cite{berard83:_scalar} that there is a metric of positive scalar curvature on \(M\).
It follows from the proof of that theorem that this metric can be chosen to be \(T^{\mathfrak{F}(P)}/T'\)-invariant.

\bibliography{pos_scalar}{}
\bibliographystyle{amsplain}
\end{document}